\documentclass[3p,12pt]{elsarticle}

\usepackage{fleqn}

\usepackage{geometry}
\usepackage{natbib}
\usepackage{pifont}
\usepackage{pifont}
\usepackage{natbib}
\usepackage{amssymb}

\usepackage{amsmath}
\newtheorem{theorem}{Theorem}[section]
\newtheorem{lemma}[theorem]{Lemma}
\newtheorem{proposition}[theorem]{Proposition}
\newtheorem{corollary}[theorem]{Corollary}

\numberwithin{equation}{section}

\journal{... }

\begin{document}

\begin{frontmatter}



\title{Approximating fixed points of asymptotically nonexpansive
mappings in Banach spaces by metric projections}


\author[h]{Hossein Dehghan}
\ead{h$_{-}$dehghan@iasbs.ac.ir, hossein.dehgan@gmail.com}

\address[h]{Department of Mathematics, Institute for
Advanced Studies in Basic Sciences (IASBS), Gava Zang, Zanjan 45137-66731, Iran
 }


\begin{abstract}
 \par In this paper, a strong convergence theorem for
asymptotically nonexpansive mappings in a uniformly convex and
smooth Banach space is proved by using metric projections. This
theorem extends and improves the recent strong convergence theorem due
to Matsushita and Takahashi [
 Appl. Math. Comput. 196 (2008) 422-425]
which was established for nonexpansive mappings.
\end{abstract}


\end{frontmatter}


\section{Introduction}
\label{}

Let $C$ be a closed convex subset of a real Banach space
$E$.  A mapping $T:C \rightarrow E$ is called nonexpansive if
$\|Tx-Ty\|\leq \|x-y\|$ for all $x,y\in C$. Also  a mapping $T : C
\rightarrow C$ is called asymptotically nonexpansive if there exists
a sequence $\{k_n\}\subset [1,\infty )$ with $k_n\to 1$ as
$n\to\infty$ such that
\begin{eqnarray*}\|T^nx-T^ny\|\leq k_n\|x-y\|\end{eqnarray*}
 for all $x,y\in C$ and
each $n\geq1$. The class of asymptotically nonexpansive mappings was
introduced by Goebel and Kirk \cite{GK} as an
important generalization of nonexpansive mappings. It was proved in
\cite{GK} that if $C$ is a nonempty bounded closed convex subset of
a real uniformly convex Banach space and $T$ is an asymptotically
nonexpansive self mapping on $C$, then $F(T)$ is nonempty closed convex subset of $C$, where $F(T)$ denotes the set of all fixed points of $T$. Strong
convergence theorems for asymptotically nonexpansive mappings have
been investigated with implicit and explicit iterative schemes (see
\cite{CZG,Sc,DGA,XA} and references therein). On the other hand, using the
metric projection, Nakajo and Takahashi \cite{NT1} introduced the
following iterative algorithm for the nonexpansive mapping $T$ in
the framework of Hilbert spaces: $x_0=x\in C$ and
\begin{eqnarray}\label{nt}
\left\{\begin{array}{l}
              y_n=\alpha_nx_n+(1-\alpha_n)Tx_n,\\
               C_n=\{z\in C : \|z-y_n\|\leq \|z-x_n\|\},\\
           Q_n=\{z\in C : \langle x_n-z,x-x_n \rangle \geq 0 \},\\
           x_{n+1}= P_{C_n\cap Q_n}x, \ \ \ \ n=0,1,2,\ldots,
         \end{array}
\right. \end{eqnarray}
  where $\{\alpha_n\}\subset [0,\alpha]$, $\alpha\in [0,1)$
 and $P_{C_n\cap Q_n}$ is the metric projection from a
Hilbert space $H$ onto $C_n \cap Q_n$. They proved that $\{x_n\}$
generated by  (\ref{nt}) converges strongly to a fixed point of $T$.
Xu \cite{Xu} extended Nakajo and Takahashi's theorem to Banach
spaces by using the  generalized projection. \par Matsushita and
Takahashi \cite{MT} recently introduced the following iterative
algorithm in the framework of Banach spaces: $x_0=x\in C$ and
\begin{eqnarray}\label{mt}
\left\{\begin{array}{l}
           C_n=\overline{co}\{z\in C : \|z-T^{}z\|\leq t_{n}\|x_{n}-T^{}x_{n}\|\},\\
           D_n=\{z\in C : \langle x_n-z,J(x-x_n) \rangle \geq 0 \},\\
           x_{n+1}= P_{ C_n\cap D_n}x, \ \ \ \ n=0,1,2,\ldots,
         \end{array}
\right.
\end{eqnarray}
where $\overline{co}D$ denotes the convex closure of the set $D$,
$J$ is normalized duality mapping, $\{t_n\}$ is a sequence in $(0,
1)$ with $t_n\to 0$, and $P_{C_n\cap D_n}$ is the metric projection
from $E$ onto $C_n\cap D_n$. Then, they proved that $\{x_n\}$
generated by (\ref{mt}) converges strongly to a fixed point of nonexpansive mapping T.
\par In this paper, motivated by (\ref{nt}) and (\ref{mt}), we
introduce the following iterative algorithm for finding fixed points
of asymptotically nonexpansive mapping $T$ in a uniformly convex and
smooth Banach space:  $x_1= x \in C$,  $C_0=D_0=C$ and
\begin{eqnarray}\label{h}
\left\{\begin{array}{l}
           C_n=\overline{co}\{z\in C_{n-1} : \|z-T^{n}z\|\leq t_{n}\|x_{n}-T^{n}x_{n}\|\},\\
           D_n=\{z\in D_{n-1} : \langle x_n-z,J(x-x_n) \rangle \geq 0 \},  \\
           x_{n+1}= P_{C_n\cap D_n}x, \ \ \ \ n=1,2,\ldots,
         \end{array}
\right.
\end{eqnarray}
where $\overline{co}D$ denotes the convex closure of the set $D$,
$J$ is normalized duality mapping, $\{t_n\}$ is a sequence in $(0,
1)$ with $t_n\to 0$, and $P_{C_n\cap D_n}$ is the metric projection
from $E$ onto $C_n\cap D_n$.\par
 The purpose of this paper is to
establish a strong convergence theorem of the iterative algorithm
(\ref{h}) for asymptotically nonexpansive mappings in a uniformly
convex and smooth Banach space.
 The results presented in this paper extend and improve the
 corresponding ones announced by Matsushita and Takahashi \cite{MT} and many others.
  \section{\textbf{Preliminaries}}
\bigskip
In this section, we recall the well-known concepts and results which
are needed to prove our main convergence  theorem. Throughout this
paper we denote by $\mathbb{N}$ the set of all positive integers.
Let $E$ be a real Banach space and let $E^*$ be the dual of $E$. We
denote the value of $x^* \in E^*$ at $x \in E$ by
$\left<x,x^*\right>$. When $\{x_n\}$ is a sequence in E, we denote
strong convergence of $\{x_n\}$ to $x \in E$ by $x_n \rightarrow x$
and weak convergence by $x_n \rightharpoonup x$. The normalized
duality mapping $J$ from $E$ to $2^{E^*}$ is defined by
\begin{eqnarray}\hspace{0cm}\nonumber J(x)=\{x^* \in E^* :
\left<x,x^*\right>=\|x\|^2=\|x^*\|^2 \}\end{eqnarray}
 for all $x \in E$. Some properties of
 duality mapping have been given in \cite{Ci}.\par
 A Banach space $E$
is said to be \emph{strictly convex} if $\|(x+y)/2\|<1$  for
all $x, y \in E$ with $\|x\|=\|y\|=1$ and $x\neq y$. A Banach space
$E$ is also said to be \emph{uniformly convex} if
$\lim_{n\to\infty}\|x_n -y_n\|= 0$ for any two sequences $\{x_n\}$
and $\{y_n\}$ in $E$ such that $\|x_n\|=\|y_n\|=1$  and
$\lim_{n\to\infty}\|x_n +y_n\|= 2$. We also know that if $E$ is a
uniformly convex Banach space, then $x_n \rightharpoonup x$ and
$\|x_n\| \rightarrow \|x\|$ imply $x_n \rightarrow x$. Let $U=\{x\in
E : \|x\|=1\}$ be the unit sphere of $E$. Then the Banach space E is
said to be \emph{smooth} if
\begin{eqnarray}\hspace{0cm}\nonumber\lim_{t\to 0}\frac{\|x+ty\|-\|x\|}{t}\end{eqnarray}
 exists for each $x, y \in U$. It is known that a Banach space $E$ is smooth if and only if
 the normalized duality mapping $J$ is single-valued. Let $C$ be a
closed convex subset of a reflexive, strictly convex and smooth
Banach space $E$. Then for any $x\in E$ there exists a unique point
$x_0\in C$ such that
$\hspace{0cm}\nonumber\|x_0-x\|= \min_{ y\in C} \|y- x\|$.
 The
mapping $P_C : E \rightarrow C$ defined by $P_Cx= x_0$ is called the
\emph{metric projection} from $E$ onto $C$. Let $x\in E$ and $u\in
C$. Then, it is known that $u=P_Cx$ if and only if
\begin{eqnarray}\label{Pro}\hspace{0cm}\left<u-y,J(x-u)\right> \geq0\end{eqnarray}
 for all $y\in C$ (see
\cite{Ta1,Ta2}).  The following proposition was proved by Bruck \cite{Br}.

\begin{proposition}
 Let $C$ be a bounded closed convex subset of a
uniformly convex Banach space $E$. Then there exists a strictly
increasing convex continuous function $\gamma :
[0,\infty)\rightarrow [0,\infty)$ with $\gamma(0)=0$ depending only
on the diameter of  $C$ such that
\begin{eqnarray}\hspace{0cm}\nonumber\gamma\left(\left\|T\left(\sum_{i=1}^n\lambda_ix_i\right)-\sum_{i=1}^n\lambda_i
Tx_i\right\|\right)\leq \max_{1\leq i<j\leq
n}(\|x_i-x_j\|-\|Tx_i-Tx_j\|)\end{eqnarray}
 holds for any
nonexpansive mapping $T : C  \rightarrow E$, any elements $x_1,\ldots,
x_n$ in $C$ and any numbers $\lambda_1,\ldots,\lambda_n\geq0$ with
$\lambda_1+\ldots+\lambda_n=1$. (Note that $\gamma$ does not depend on
$T$.)
\end{proposition}
 \begin{corollary} Under the same suppositions as in Proposition 2.1, there exists a strictly
increasing convex continuous function $\gamma :
[0,\infty)\rightarrow [0,\infty)$ with $\gamma(0)=0$ such that
\begin{eqnarray*}\hspace{0cm}\nonumber
\gamma\left(\frac{1}{k_m}\left\|T^m\left(\sum_{i=1}^n\lambda_ix_i\right)-\sum_{i=1}^n
\lambda_iT^mx_i\right\|\right)\leq \max_{1\leq i<j\leq
n}\left(\|x_i-x_j\|-\frac{1}{k_m}\left\|T^mx_i-T^mx_j\right\|\right)
\end{eqnarray*}
   for any asymptotically nonexpansive mapping
$T: C\rightarrow C$ with $\{k_n\}$,
  any elements $x_1,\ldots,
x_n$ in $C$, any numbers $\lambda_1,\ldots,\lambda_n\geq0$ with
$\lambda_1+\cdots+\lambda_n=1$ and each  $m\geq1$.
\end{corollary}
\textbf{Proof.} Define the mapping $S_m : C \rightarrow E$ as $S_m
x=1/k_mT^mx,$ for all $x\in C$ and each $m\geq1$. Then $S_m$ is
nonexpansive for all $m\geq1$. From proposition 2.1, there exists a
strictly increasing convex continuous function $\gamma :
[0,\infty)\rightarrow [0,\infty)$ with $\gamma(0)=0$ such that
\begin{eqnarray*}\hspace{0cm}\nonumber\gamma\left(\left\|S_m\left(\sum_{j=1}^n\lambda_jx_j\right)-
\sum_{j=1}^n\lambda_jS_mx_j\right\|\right) \leq \max_{1\leq j<k\leq
n}(\|x_j-x_k\|-\|S_mx_j-S_mx_k\|)\end{eqnarray*}
 for all
$m\geq1$. Thus, by using the definition of $S_m$, we obtain the
desired conclusion. $\Box$
\begin{lemma}\cite[Lemma 1.6]{CZG}
 Let $E$ be a uniformly convex Banach space, $C$ be a nonempty closed
convex subset of $E$ and $T:C\rightarrow C$ be an asymptotically
nonexpansive mapping. Then $(I-T)$ is demiclosed at $0$, i.e., if
$x_n\rightharpoonup x$  and $x_n-Tx_n\to 0$, then $x\in F(T)$, where
$F(T)$ is the set of all fixed points of $T$.
\end{lemma}
\section{\textbf{Strong convergence theorem}}
 In this section, we study the iterative
algorithm (\ref{h}) for finding fixed points of asymptotically
nonexpansive mappings in a uniformly convex and smooth Banach space.
 We first prove that the sequence $\{x_n\}$ generated by (\ref{h}) is
well-defined. Then, we prove that $\{x_n\}$ converges strongly to
$P_{F(T)}x$, where $P_{F(T)}$ is the metric projection from $E$ onto
$F(T)$.
\begin{lemma}
Let $C$ be a nonempty closed convex subset of a reflexive, strictly
convex and smooth Banach space $E$ and let $T : C \rightarrow C$ be
an asymptotically  nonexpansive mapping. If $F(T)\neq \emptyset$,
then the sequence $\{x_n\}$  generated by (\ref{h}) is well-defined.
\end{lemma}
\textbf{Proof.} It is easy to check that $C_n \cap D_n$ is closed
and convex and $F(T)\subset C_n$ for each $n \in\mathbb{ N}$.
Moreover $D_1=C$ and so
 $F(T)\subset
C_1\cap D_1$. Suppose $F(T)\subset C_k \cap D_k $ for $k \in\mathbb{
N}$. Then, there exists a unique element $x_{k+1}þ\in C_k \cap D_k$
such that $x_{k+1}=P_{C_k\cap D_k} x$. If $u \in F(T)$, then it
follows from (2.1) that
\begin{eqnarray}\hspace{0cm}\nonumber\left<x_{k+1}-u,J(x-x_{k+1})\right>\geq0,\end{eqnarray}
which implies $u\in  D_{k+1}$. Therefore  $F(T) \subset C_{k+1} \cap
D_{k+1} $.
 By mathematical
induction, we obtain that $F(T) \subset C_n \cap D_n $ for all $n
\in \mathbb{N}$. Therefore, $\{x_n\}$ is well-defined. $\Box$\\
\par In order to prove our main result, the following lemma is needed.
\begin{lemma}
 Let $C$ be a nonempty bounded closed convex subset of a uniformly convex and
smooth Banach space $E$ and let $T : C \rightarrow C$ be an
asymptotically nonexpansive mapping with $\{k_n\}$ and  $\{x_n\}$ be
the sequence generated by (\ref{h}). Then for any $k\in\mathbb{N}$,
\begin{eqnarray}\hspace{0cm}\nonumber\lim_{n\to\infty}\left\|x_{n}-T^{n-k}x_{n}\right\|=0.
\end{eqnarray}
\end{lemma}
\textbf{Proof.}
 Fix $k\in\mathbb{N}$ and put $m=n-k$. Since  $x_{nþ} =P_{
C_{n-1}\cap D_{n-1}}x$, we have $x_{nþ} \in C_{n-1}\subseteq\cdots\subseteq
C_{m}$.
 Since $t_{m} > 0$, there exist
 $y_1,\ldots,y_N\in
C$ and $\lambda_1,\ldots,\lambda_N\geq0$ with
$\lambda_1+\cdots+\lambda_N=1$ such that
\begin{eqnarray}\label{1}\hspace{-0.5cm}\left\|x_{n}-\sum_{i=1}^N\lambda_iy_i\right\|<t_{m},\end{eqnarray}
 and $\|y_i+T^{m}y_i\|\leq t_{m}\|x_{m}-T^{m}x_{m}\|$
for all $i\in \{1,\ldots,N\}$. Put  $M= \sup_{x \in C} \|x\| $,
$u=P_{F(T)}x$ and $r_0=\sup_{n\geq1}(1+k_n)\|x_n-u\|$. Since $C$ and
$\{k_n\}$ are bounded,  (\ref{1}) implies
\begin{eqnarray}\label{3}\hspace{-.5cm}\left\|x_{n}-\frac{1}{k_{m}}\sum_{i=1}^N\lambda_iy_i\right\|
\leq \left(1-\frac{1}{k_{m}}\right)\|x_{n}\|
+\frac{1}{k_{m}}\left\|x_{n}-\sum_{i=1}^N\lambda_iy_i\right\|\leq\left(1-\frac{1}{k_{m}}\right)M+t_{m},\end{eqnarray} and $\left\|y_i-T^{m}y_i\right\|\leq
t_{m}\|x_{m}-T^{m}x_{m}\|\leq t_{m}(1+k_{m})\|x_{m}-u\|\leq
r_0t_{m}$ for all $i\in \{1,\ldots,N\}$. Therefore
\begin{eqnarray}\hspace{-0.4cm}
\label{4}\left\|y_i-\frac{1}{k_{m}}T^{m}y_i\right\|\leq
\left(1-\frac{1} {k_{m}}\right)M+r_0t_{m}\end{eqnarray} for all
$i\in \{1,\ldots,N\}$. Moreover, asymptotically nonexpansiveness of $T$ and (\ref{1}) give that
\begin{eqnarray}\hspace{-0.4cm}
\label{44}\left\|\frac{1}{k_{m}}T^{m}\left(\sum_{i=1}^N\lambda_iy_i\right)
-T^{m}x_{n}\right\|\leq
\left(1-\frac{1} {k_{m}}\right)M+t_{m}.\end{eqnarray}
It follows from Corollary 2.2,  (\ref{3}), (\ref{4}) and
(\ref{44}) that
\begin{eqnarray}\label{5}\hspace{-0.4cm}
\nonumber\left\|x_{n}-T^{m}x_{n}\right\|\hspace{0cm}&\leq&\hspace{0cm}\left\|
x_{n}-\frac{1}{k_{m}}\sum_{i=1}^N\lambda_iy_i\right\|
+\frac{1}{k_{m}}\left\|\sum_{i=1}^N\lambda_i\left(y_i-T^{m}y_i\right)\right\|\\
\nonumber&\ \ \ &+\frac{1}{k_{m}}\left\|\sum_{i=1}^N\lambda_iT^{m}y_i-T^{m}\left(
\sum_{i=1}^N\lambda_iy_i \right)\right\|\\
\nonumber&\ \ \ &+\left\|\frac{1}{k_{m}}T^{m}\left(\sum_{i=1}^N\lambda_iy_i\right)
-T^{m}x_{n}\right\|\\
\nonumber&\leq&
2\left(1-\frac{1}{k_{m}}\right)M+2t_{m}+\frac{r_0t_{m}}{k_m}\\
\nonumber&\ \ \ &+\gamma_{}^{-1}\left(\max_{1\leq i<j\leq
N}\left(\|y_i-y_j\|-\frac{1}{k_{m}}\left\|T^{m}y_i-T^{m}y_j\right\|\right)\right)\\
\nonumber&
\leq&2\left(1-\frac{1}{k_{m}}\right)M+2t_{m}+\frac{r_0t_{m}}{k_m}\\
\nonumber&\ \ \ &+\gamma_{}^{-1}\left(\max_{1\leq i<j\leq
N}\left(\left\|y_i-\frac{1}{k_{m}}T^{m}y_i\right\|+\left\|y_j-\frac{1}{k_{m}}T^{m}y_j
\right\|\right)\right)\\
\nonumber& \leq&
2\left(1-\frac{1}{k_{m}}\right)M+2t_{m}+\frac{r_0t_{m}}{k_m}
+\gamma_{}^{-1}\left(2\left(
1-\frac{1}{k_{m}} \right)M+2r_0t_{m}\right).
\end{eqnarray}
Since $\lim_{n\to\infty}k_n=1$ and $\lim_{n\to\infty}t_n=0$, it
follows from  the last inequality that
$\lim_{n\to
\infty}\|x_{n}-T^mx_{n}\|= 0$. This completes the
proof.  $\Box$
\begin{theorem} Let $C$ be a nonempty bounded closed convex subset of a
uniformly convex and smooth Banach space $E$ and let $T : C
\rightarrow C$ be an asymptotically nonexpansive mapping and let
$\{x_n\}$ be the sequence generated by (\ref{h}). Then $\{x_n\}$
converges strongly to the element $P_{F(T)}x$ of $F(T)$, where
$P_{F(T)}$ is the metric projection from $E$ onto $F(T)$.
\end{theorem}
\textbf{Proof.}
 Put $u=P_{F(T)}x$. Since $F(T)\subset C_{n}\cap D_{n}$ and $x_{n+1}= P_{C_{n}\cap D_{n}}x$,
we have that \begin{eqnarray}\label{6}\hspace{0cm} \|x-
x_{n+1}\|\leq\|x- u\|\end{eqnarray}
 for all $n\in \mathbb{ N}$. By Lemma 3.2, we have
\begin{eqnarray*}\hspace{0cm}\left\|x_{n}-Tx_{n}\right\|\hspace{-.2cm}&\leq&\hspace{-.2cm}\left\|x_{n}-T^{n-1}x_{n}
\right\|+\left\|T^{n-1}x_{n}-Tx_{n}\right\|\\
&\leq&\hspace{-.2cm}
\left\|x_{n}-T^{n-1}x_{n}\right\|+k_1\left\|T^{n-2}x_{n}-x_{n}
\right\|\rightarrow 0 \ \ \ as\ n\to\infty.\end{eqnarray*} Since
$\{x_n\}$ is bounded, there exists $\{x_{n_i}\}\subset \{x_n\}$ such
that $x_{n_i}\rightharpoonup v$. It follows from Lemma 2.3 that
$v\in F(T)$. From the weakly lower semicontinuity of norm and
(\ref{6}),
 we obtain
\begin{eqnarray*}\hspace{0cm}\|x-u\|\leq\|x-v\|\leq \liminf_{i\to\infty}\|x-x_{n_i}\|\leq
  \limsup_{i\to\infty}\|x-x_{n_i}\|\leq\|x-u\|.\end{eqnarray*}
This together with the uniqueness of $P_{F(T)}x$, implies $u=v$, and
hence $x_{n_i}\rightharpoonup u $.
 Therefore, we obtain $x_n\rightharpoonup  u$. Furthermore, we have that
\begin{eqnarray*}\hspace{0cm}\lim_{ n\to\infty} \|x - x_n\|=\| x -u\|.\end{eqnarray*}
 Since E is uniformly convex, we have
that $x- x_n \to x -u$. It follows that $x_n \to u$. This completes
the proof. $\Box$\\\\











\end{document}